\documentclass[a4paper,11pt]{amsart}
\usepackage{amssymb,amsmath,amsthm,graphicx}
\usepackage{enumerate}

\newtheorem{theorem}{Theorem}[section]
\newtheorem*{theorem*}{Theorem}
\newtheorem{lemma}[theorem]{Lemma}
\newtheorem{proposition}[theorem]{Proposition}
\newtheorem{corollary}[theorem]{Corollary}

\theoremstyle{remark}

\theoremstyle{definition}
\newtheorem{example}[theorem]{Example}
\newtheorem{definition}[theorem]{Definition}

\begin{document}
\title[Metric properties in Berggren tree of PPT]{Metric properties in Berggren tree of primitive Pythagorean triples}

\author[L.~Jani\v{c}kov\'a]{Lucia Jani\v{c}kov\'a}

\author[E. Cs\'ok\'asi]{Evelin Cs\'ok\'asi}

\begin{abstract}
A Pythagorean triple is a triple of positive integers $(x,y,z)$ such that $x^2+y^2=z^2$. If $x,y$ are coprime and $x$ is odd, then it is called a primitive Pythagorean triple.  Berggren showed that every primitive Pythagorean triple can be  generated from triple $(3,4,5)$ using multiplication by uniquely number and order of three $3\times3$ matrices, which yields a ternary tree of triplets. In this paper, we present some metric properties of triples in Berggren tree.  Firstly, we consider  primitive Pythagorean triple as lengths of sides of right triangles and secondly, we consider them as coordinates of points in three-dimensional space.
\end{abstract}
\maketitle


\section{Introduction}

There are various ways of generating integer solutions of Pythagorean equation. Euclid's formula, Berggren tree and Price tree are among the best known ones, however there are many others. Some examples can be found e.g. in \cite{Lbc}, \cite{Emgr}. The matrices used to generate Berggren tree produce coprime solutions - the primitive Pythagorean triples, which can be viewed as lengths of the sides of Pythagorean triangles.

Some properties of the Pythagorean triangles were already described. E.g., the inradius \cite{2006}, triples with common lengths of leg \cite{legs} or height of primitive Pythagorean triples (the difference between length of hypotenuse and length of even leg)  \cite{height}. 

In this paper, we present our results concerning metric properties of triangles with lengths of sides corresponding to primitive Pythagorean triples and of triangles formed by descendants in Berggren tree. We focus on some specific sequences of primitive Pythagorean triples in Berggren tree and their inradii and radii of their circumcircles.

Further, the primitive Pythagorean triples can be viewed as coordinates of the points in the 3-dimensional real space.  We explore some of the properties of the points which coordinates are primitive Pythagorean triples.

 
 \section{Preliminary}
First, we present some basic notations and definitions which we use through the paper. 

Let us denote $\mathbb{N}:= \{1, 2, 3, \dots \}$ and $\mathbb{N}_0:= \mathbb{N} \cup \{0\}$.  

For a matrix $M$, denote its determinant as $\det(M)$ and its transpose as $M^\top.$

Let $E$ denote the unit matrix of the size $3$.

For a vector $\vec{u}$, we denote its size as $|\vec{u}|$.

\begin{definition}
A triple of positive integers $(x,y,z)$ is called \textit{a Pythagorean triple} if $x^2+y^2=z^2.$ Moreover, if $\gcd(x,y)=1$ and $x$ is odd, then we call $(x,y,z)$ \textit{a primitive Pythagorean triple}.
\end{definition}

If $(x,y,z)$ is a Pythagorean triple, we say that the triangle \textit{corresponds to} this triple if its sides have lengths $x,y,z$. A triangle corresponding to a Pythagorean triple is called a \textit{Pythagorean triangle.}

According to Berggren, every primitive Pythagorean triple $(x,y,z)$ with $y$ even can
be generated from the triple $(3,4,5)$ by unique 3-fold ascent using the three matrices $A, B, C$  \cite{hall}:
$$
A= \begin{pmatrix}
1& -2& 2\\
2&-1& 2\\
2&-2& 3
\end{pmatrix}, \quad
B= \begin{pmatrix}
1&  2& 2\\
2& 1& 2\\
2& 2& 3
\end{pmatrix},\quad 
C= \begin{pmatrix}
-1&  2& 2\\
-2& 1& 2\\
-2& 2& 3
\end{pmatrix}.
$$
Then for every $n\in\mathbb{N}$:
\begin{align*}
A^n&= \begin{pmatrix}
1& -2n& 2n\\
2n& 1-2n^2& 2n^2\\
2n&-2n^2& 2n^2+1
\end{pmatrix},\\
B^n&= \begin{pmatrix}
\frac{(-1)^n}{2}+b_1&  \frac{(-1)^{n+1}}{2}+ b_1& b_2\\
\frac{(-1)^{n+1}}{2}+ b_1& \frac{(-1)^n}{2}+b_1& b_2\\
b_2& b_2& 2b_1 
\end{pmatrix},\\
C^n&= \begin{pmatrix}
1-2n^2&  2n& 2n^2\\
-2n& 1& 2n\\
-2n^2& 2n& 2n^2+1
\end{pmatrix},
\end{align*}
where $b_1 = \frac{1}{4}[(3-2 \sqrt{2})^n+(3+2 \sqrt{2})^n], \quad
b_2 = -\frac{\left(3-2 \sqrt{2}\right)^n}{2 \sqrt{2}}+\frac{\left(3+2 \sqrt{2}\right)^n}{2 \sqrt{2}}.$
\begin{definition}
	Let $P$ be a primitive Pythagorean triple. We say that $P$ is \textit{the parent} of the triples $AP^\top, BP^\top, CP^\top.$ The triples $AP^\top, BP^\top, CP^\top$ are called \textit{the descendants} of $P$ in Berggren tree.
\end{definition}

It is often useful to express the primitive Pythagorean triples by Euclid's formula \cite{maor}:

\begin{theorem}[Euclid's formula]
	Triple $(x,y,z)$ is a primitive Pythagorean triple with odd $x$  if and only if there exist  $m,n\in\mathbb{N}$ such that $m>n, \gcd(m,n)=1$ and $$x=m^2-n^2, \quad\quad y= 2mn, \quad\quad z= n^2+m^2.$$
	\end{theorem}

Another approach is mentioned in \cite{Emgr}: If $(x,y,z)$ is a Pythagorean triple then there exist $m,n\in \mathbb{N}$ such that 
$$
x = \frac{3-(-1)^m}{2}mn+m,\quad\quad
y = \frac{x^2-m^2}{2m},\quad\quad
z = \frac{x^2+m^2}{2m}.
$$
Then the tripple $(x,y,z)$ is denoted $F(m,n)$. 

\begin{lemma} \label{Fmn}
For every $n\in \mathbb{N}$ it holds that $A^{n-1}(3,4,5)^\top =
F(1,n)^\top.$
\end{lemma}
\begin{proof} For every $n\in \mathbb{N}$:
	\begin{align*}
	A^{n-1}(3,4,5)^\top &= \begin{pmatrix}
	1& -2(n-1)& 2(n-1)\\
	2(n-1)& 1-2(n-1)^2& 2(n-1)^2\\
	2(n-1)&-2(n-1)^2& 2(n-1)^2+1
	\end{pmatrix}\cdot\begin{pmatrix}
	3\\4\\5
	\end{pmatrix} \\
	&= (2n+1,2n^2+2n,2n^2+2n+1)^\top \\
	&= F(1,n)^\top.
	\end{align*}
	
\end{proof}
\section{Pythagorean triangles}

In \cite{Lbc}, some properties of incircle of a Pythagorean triangle were proved. In this section, we present some further results related to incircle and excircle of a Pythagorean triangle.

Clearly, if $(x,y,z)$ is a primitive Pythagorean triple, then the inradius of the corresponding triangle is $\frac{x+y-z}{2}.$

\begin{proposition}
Let $(x,y,z)$ be a primitive Pythagorean triple and let $r$ be the inradius of the corresponding triangle. Then the inradius of the triangle corresponding to triple $A\cdot(x,y,z)^{\top}, B\cdot(x,y,z)^{\top}$ and $C\cdot(x,y,z)^{\top}$ is \begin{align*}
r_A&=r-y+z,\\ r_B&=r+z, \\r_C&=r-x+z,
\end{align*}respectively.
\end{proposition}
For proof, see \cite{Lbc}.

\begin{proposition}
If $n\in \mathbb{N}$ then the inradius of the triangle corresponding to triple $A^n\cdot(3,4,5)^{\top}$ is $r_n=n+1.$
\end{proposition}

\begin{proof}
It is easy to show that $A^{n}\cdot(3,4,5)^{\top} = (2n+3,2n^2+6n+4,2n^2+6n+5)^{\top},$ hence $r_n=\frac{(2n+3)+(2n^2+6n+4)-(2n^2+6n+5)}{2}=n+1.$
\end{proof}

\begin{corollary}
	For every $r\in \mathbb{N},$ there exists a primitive Pythagorean triple such that the inradius of the triangle corresponding to this triple is $r$.
\end{corollary}

Notice that the inradius does not determine the primitive Pythagorean triple unambiguously. For example, the triangles corresponding to $(21,20,29)$ and $(13,84,85)$ have the same inradius $r=6.$ This leads to natural question - how many primitive Pythagorean triples have the same inradius? In 2006, Robbins answered this question in \cite{2006}.

\begin{theorem} \label{polomer}
Let $\omega(r)$ be the number of prime divisors of the integer $r$ and let $I(r)$ be the number of primitive Pythagorean triples with the inradius $r.$ Then  

\begin{align*}
I(r) &= \left\{\begin{array}{lll}
2^{\omega(r)} \quad &\text{if}~r~\text{is odd}, \\
2^{\omega(r)-1} \quad &\text{if}~r~\text{is even}.
\end{array}\right.\\
\end{align*} 
\end{theorem}

In 2016, Omland presented a simplified proof, see \cite{inradius}.  It suggests an easy way of finding all primitive Pythagorean triples with given inradius. The following example illustrates the case for inradius 35.
 
\begin{example}
Let $r=35$. Its prime decomposition is $35=5^1\cdot 7^1$, hence $\omega(r)=2$ According to Theorem $\ref{polomer},$ there are $2^{\omega(r)}=4$ primitive Pythagorean triples with inradius 35. It is easy to see:

\begin{figure}[h!t]
	\centering
	\includegraphics[width=0.8\textwidth]{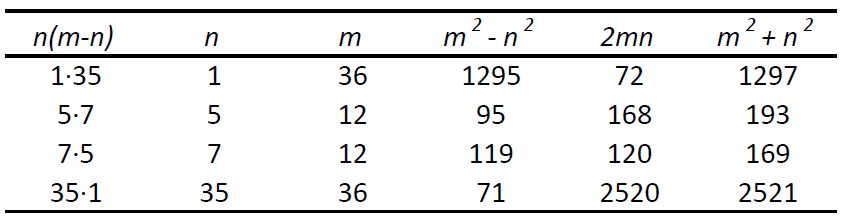}
	\label{example1}
	\caption{Pythagorean triples with inradius 35}
\end{figure}

Hence, $(1295,72,1297), (95,168,193), (119,120,169), (71,2520,2521)$ are all  primitive Pythagorean triples with inradius 35.
\end{example}

\begin{proposition} \label{bn}
If $n\in \mathbb{N}$ then the inradius of the triangle corresponding to triple $B^n\cdot(3,4,5)^{\top}$ is $r_n=\frac{(3+2\sqrt{2})^{n+1}-(3-2\sqrt{2})^{n+1}}{4\sqrt{2}}.$
\end{proposition}

\begin{proof}
Analogously like above, $$B^{n}\cdot(3,4,5)^{\top}  
= \begin{pmatrix}
-\frac12(-1)^n+\frac{7\sqrt{2}-10}{4\sqrt{2}}(3-2\sqrt{2})^n+ \frac{7\sqrt{2}+10}{4\sqrt{2}}(3+2\sqrt{2})^n \\
\frac12(-1)^n+\frac{7\sqrt{2}-10}{4\sqrt{2}}(3-2\sqrt{2})^n+ \frac{7\sqrt{2}+10}{4\sqrt{2}}(3+2\sqrt{2})^n \\
\frac{(5\sqrt{2}-7)(3-2\sqrt{2})^n+(5\sqrt{2}+7)(3+2\sqrt{2})^n}{2\sqrt{2}}
\end{pmatrix},\\ $$ hence
 
\begin{align*}
 r_n&=\frac{\frac{7\sqrt{2}-10}{2\sqrt{2}}(3-2\sqrt{2})^n+ \frac{7\sqrt{2}+10}{2\sqrt{2}}(3+2\sqrt{2})^n-
 \frac{(5\sqrt{2}-7)(3-2\sqrt{2})^n+(5\sqrt{2}+7)(3+2\sqrt{2})^n}{2\sqrt{2}}}{2}
 \\&=\frac{(2\sqrt{2}-3)(3-2\sqrt{2})^{n}+(2\sqrt{2}+3)(3+2\sqrt{2})^{n}}{4\sqrt{2}}
 \\&=\frac{(3+2\sqrt{2})^{n+1}-(3-2\sqrt{2})^{n+1}}{4\sqrt{2}}.
\end{align*}
\end{proof}

\begin{proposition}
If $n\in \mathbb{N}$ then the inradius of the triangle corresponding to triple $C^n\cdot(3,4,5)^{\top}$ is $r_n=2n+1.$
\end{proposition}

\begin{proof}
Analogously like above, $C^{n}\cdot (3,4,5)^{\top} = (4n^2+8n+3, 4n+4, 4n^2+8n+5)^{\top},$ hence $r_n=\frac{(n^2+8n+3)+(4n+4)-(4n^2+8n+5)}{2}=2n+1.$
\end{proof}

Similarly, it is easy to show that if $(x,y,z)$ is a primitive Pythagorean triple, then the radius of the circumcircle of the corresponding triangle is $\frac{z}{2}.$ 

\begin{proposition}
Let $(x,y,z)$ be a primitive Pythagorean triple and let $R$ be the radius of the circumcircle of the corresponding triangle. Then the radius of the circumcircle of the  triangle corresponding to triple $A\cdot(x,y,z)^{\top}, B\cdot(x,y,z)^{\top}$ and $C\cdot(x,y,z)^{\top}$ is \begin{align*}
R_A&=x-y+3R,\\ R_B&=x+y+3R, \\R_C&=-x+y+3R,
\end{align*}respectively.
\end{proposition}

\begin{proof}
Follows directly from $R= \frac{z}{2}.$
\end{proof}

\begin{proposition}
If $n\in \mathbb{N}$ then the radius of the circumcircle of the triangle corresponding to triple $A^n\cdot(3,4,5)^{\top}$ is $R_n=n^2+3n+\frac52.$
\end{proposition}

\begin{proof}
From $A^{n}\cdot (3,4,5)^{\top} = (2n+3,2n^2+6n+4,2n^2+6n+5)^{\top},$ it follows that $R_n=\frac{2n^2+6n+5}{2}=n^2+3n+\frac52.$
\end{proof}

\begin{proposition} 
If $n\in \mathbb{N}$ then the radius of the circumcircle of the triangle corresponding to triple $B^n\cdot(3,4,5)^{\top}$ is $ R_n=
 \frac{(5\sqrt{2}-7)(3-2\sqrt{2})^n+(5\sqrt{2}+7)(3+2\sqrt{2})^n}{4\sqrt{2}}.$
\end{proposition}

\begin{proof}
Analogously like in the proof of Proposition \ref{bn}, we compute $B^{n}\cdot (3,4,5)^{\top}$ and then $$
 R_n=\frac12 \cdot 
 \frac{(5\sqrt{2}-7)(3-2\sqrt{2})^n+(5\sqrt{2}+7)(3+2\sqrt{2})^n}{2\sqrt{2}}.
$$
\end{proof}

\begin{proposition}
If $n\in \mathbb{N}$ then the radius of circumcircle of the triangle corresponding to triple $C^n\cdot(3,4,5)^{\top}$ is $R_n=2n^2+4n+\frac52.$
\end{proposition}

\begin{proof}
From $C^{n}\cdot (3,4,5)^{\top} = (4n^2+8n+3, 4n+4, 4n^2+8n+5)^{\top},$ it follows that $R_n=\frac{4n^2+8n+5}{2}=2n^2+4n+\frac52.$
\end{proof}


 \section{Relations in the triangle of descendant triples}

The primitive Pythagorean triples can be viewed as points in the 3-dimensional Euclidean space. In this section, we show that the  descendants of any primitive Pythagorean triple in Berggren tree are vertices of a triangle, and we present our results related to the geometric relations in these triangles. Namely, we study the triangle defined by the triple of descendants.

\begin{proposition} \label{rovina}
Let $P$ be a primitive Pythagorean triple. Then the points with the coordinates $AP^\top, BP^\top, CP^\top$ are not collinear.
\end{proposition}

\begin{proof}
Let $P=(x,y,z)$. We consider the vectors $\vec{u}=BP^\top - AP^\top$ and $\vec{v}=CP^\top - AP^\top$. It is easy to show that 
\begin{align*}
\vec{u}&= (B-A)P^\top = \begin{pmatrix}
0 & 4 & 0  \\
0 & 2 & 0  \\
0 & 4 & 0
\end{pmatrix}\cdot \begin{pmatrix}
x  \\
y  \\
z
\end{pmatrix} =
\begin{pmatrix}
4y  \\
2y \\
4y
\end{pmatrix},
\\ 
\vec{v}&= (C-A)P^\top = \begin{pmatrix}
-2 & 4 & 0  \\
-4 & 2 & 0  \\
-4 & 4 & 0
\end{pmatrix}\cdot \begin{pmatrix}
x  \\
y  \\
z
\end{pmatrix} =
\begin{pmatrix}
-2x+4y  \\
-4x+2y \\
-4x+4y
\end{pmatrix}.
\end{align*}

By way of contradiction, assume that these vectors are linearly dependent, i.e., that there exists nonzero real number $k$  such that $\vec{u}=k\vec{v}.$ Then
 $(4y,2y,4y)=k(-2x+4y,-4x+2y,-4x+4y)$, hence $4y=k(-2x+4y)$ and $4y=k(-4x+4y)$. This implies that 
 \begin{align*}
 k(-2x+4y)&=k(-4x+4y) \\
 -2x+4y&=-4x+4y \\
 x&=2x
 \end{align*}
which is a contradiction. Therefore, $\vec{u}, \vec{v}$ are linearly independent and the points $AP^\top, BP^\top, CP^\top$ are not collinear.
\end{proof}

This implies that the points with the coordinates $AP^\top, BP^\top, CP^\top$ determine a plane (a triangle). In the following propositions, we express the equation of this plane (the area of this triangle).

\begin{proposition} \label{rovina2}
Let $P=(x,y,z)$ be a primitive Pythagorean triplet. Then the points with the coordinates $AP^\top, BP^\top, CP^\top$ belong to plane $2a+2b-3c+z=0.$
\end{proposition}

\begin{proof}
Analogously like in the previous proof, we consider vectors $\vec{u}=BP^\top - AP^\top$ and $\vec{v}=CP^\top - AP^\top$.  Firstly, we compute the normal vector of the wanted plane as cross product of these vectors:  $\vec{u} \times \vec{v}= (8xy, 8xy, -12xy) \thickapprox (2,2,-3).$ Hence, $2a+2b-3c+d=0$ is the equation of wanted plane for some $d\in \mathbb{R}.$ Since $AP^\top$ belongs to this plane, we get $2(x-2y+2z)+2(2x-y+2z)-3(2x-2y+3z)+d=0$ which yields $d=z.$ Therefore, the equation of the wanted plane is $2a+2b-3c+z=0.$
\end{proof}

\begin{proposition} \label{trojuholnik}
Let $P=(x,y,z)$ be a primitive Pythagorean triple. Then the points with the coordinates $AP^\top, BP^\top, CP^\top$ are vertices of a triangle with the area $2xy\sqrt{17}.$
\end{proposition}

\begin{proof}
To compute the area of the triangle $AP^\top, BP^\top, CP^\top$, we use the vectors 
 $\vec{u}=BP^\top - AP^\top$ and $\vec{v}=CP^\top - AP^\top$.  According to the previous proof,   $\vec{u} \times \vec{v}= (8xy, 8xy, -12xy)$, therefore the area of the triangle is  
 $$\frac{|\vec{u} \times \vec{v}|}{2}= \frac{\sqrt{64x^2y^2+64x^2y^2+144x^2y^2}}{2}=\frac{\sqrt{272x^2y^2}}{2}=2xy\sqrt{17}.$$
\end{proof}

Primitive pythagorean triple can be viewed as a right triangle and the points corresponding to the descendants of a PPT in Beggren tree also form a triangle. Natural question arises - can this triangle of descendants be also a right triangle? The following proposition offers the answer to this question.  

\begin{proposition} \label{4.12}
Let $P$ be a primitive Pythagorean triple. The points with coordinates $AP^\top, BP^\top, CP^\top$ form a non-right triangle.
\end{proposition}

\begin{proof}
Let $P=(x,y,z)$ and let the vectors $\vec{u}, \vec{v}, \vec{w}$ be as follows: $$\vec{u}=BP^\top-AP^\top, 
\vec{v}=CP^\top-AP^\top, 
\vec{w}=CP^\top-BP^\top.$$
To show that the triangle with the vertices $AP^\top, BP^\top, CP^\top$ is a non-right triangle, it is sufficient to show that $\vec{u}\cdot\vec{v}\not=0, \vec{u}\cdot\vec{w}\not=0, \vec{v}\cdot\vec{w}\not=0.$ According to the proof of Proposition \ref{rovina}, $\vec{u}=(4y,2y,4y)$ and $\vec{v}=(-2x+4y,-4x+2y,-4x+4y)$. 
Analogously, 
\begin{align*}
\vec{w}&=(C-B)P^\top = \begin{pmatrix}
-2 & 0 & 0  \\
-4 & 0 & 0  \\
-4 & 0 & 0
\end{pmatrix}\cdot \begin{pmatrix}
x  \\
y  \\
z
\end{pmatrix} =
\begin{pmatrix}
-2x  \\
-4x \\
-4x
\end{pmatrix}.
\end{align*}
It follows that
$\vec{u}\cdot\vec{v} = (4y,2y,4y)\cdot (-2x+4y,-4x+2z,-4x+4y) =-32xy+36y^2.$ By way of contradiction, assume that $\vec{u}\cdot\vec{v} = 0$. Since $x,y >0$, we get $y = \frac{8x}{9}.$ However, $(x,y,z)$ is a primitive Pythagorean triple which yields
\begin{align*}
x^2+\frac{64x^2}{81}&=z^2
\\
\sqrt{145} \frac{x}{9}&=z
\end{align*}
which is a contradiction with $z\in\mathbb{N}.$

Similarly, $\vec{u}\cdot\vec{w} = (4y,2y,4y)\cdot (-2x,-4x,-4x) = -32 xy$ and from $x,y >0$ it follows that $\vec{u}\cdot\vec{w}\not=0$.

Finally, $\vec{v}\cdot\vec{w} = (-2x+4y,-4x+2z,-4x+4y)\cdot (-2x,-4x,-4x) =36x^2-32xy.$ If $\vec{v}\cdot\vec{w}=0$, then we get a contradiction analogously to the case $\vec{u}\cdot\vec{v}.$ 
\end{proof}

Using the results from this section, it is easy to determine the lengths of inradius and radius of the circumcircle of of triangle whose vertics are descendants of a Pythagorean triple. 

\begin{lemma} \label{use}
Let $P$ be a primitive Pythagorean triple. Let the vectors $\vec{u}, \vec{v}, \vec{w}$ be as follows: $\vec{u}=BP^\top-AP^\top, 
\vec{v}=CP^\top-AP^\top, 
\vec{w}=CP^\top-BP^\top.$ Then the triangle with the vertices  $A P^\top, B P^\top, C P^\top$ has inradius $$r= \frac{|\vec{u}\times \vec{v}|}{|\vec{u}|+ |\vec{v}|+ |\vec{w}|}$$ 

and radius of the circumcircle 
 $$R= \frac{|\vec{u}|\cdot |\vec{v}|\cdot |\vec{w}|}{2|\vec{u}\times \vec{v}|}.$$ 
\end{lemma}

\begin{proof}
It is clear to see that $|\vec{u}|,|\vec{v}|, |\vec{w}|$ are the lengths of the sides of this triangle. We use well  known formulas for area of a triangle 
$$S = \frac{|\vec{u}\times \vec{v}|}{2}, \quad S=\frac{|\vec{u}|+ |\vec{v}|+ |\vec{w}|}{2}\cdot r, \quad S=\frac{|\vec{u}|\cdot |\vec{v}|\cdot |\vec{w}|}{4R}$$ 
to express the inradius $r$ and the radius of the circumcircle  $R$, respectively.
\end{proof}

\begin{proposition}
Let $P=(x,y,z)$ be a primitive Pythagorean triple.  Then the triangle with the vertices  $A P^\top, B P^\top, C P^\top$ has inradius 
$$r= \frac{2xy\sqrt{17}}{3x+3y+\sqrt{9x^2-16xy+9y^2}}$$ 

and radius of the circumcircle 
 $$R=\frac{9\sqrt{9x^2-16xy+9y^2}}{\sqrt{17}}.$$ 
\end{proposition}

\begin{proof}
According to the proofs of Proposition \ref{trojuholnik} and Proposition \ref{4.12}, it holds that 
\begin{align*}
 &\vec{u}=(4y,2y,4y),\\ &\vec{v}=(-2x+4y,-4x+2y,-4x+4y),\\ &\vec{w}=(-2x,-4x,-4x), \\  &|\vec{u}\times\vec{v}|=\sqrt{272x^2y^2}=4xy\sqrt{17}. 
\end{align*}
Therefore,  
\begin{align*}
 &|\vec{u}|=\sqrt{16y^2+4y^2+16y^2}=6y,\\
 &|\vec{v}|=\sqrt{36x^2-64xy+36y^2}=2\sqrt{9x^2-16xy+9y^2},\\ 
 &|\vec{w}|=\sqrt{4x^2+16x^2+16x^2}=6x. 
\end{align*}
Then according to the Lemma \ref{use}, the inradius is $$r=\frac{4xy\sqrt{17}}{6x+6y+2\sqrt{9x^2-16xy+9y^2}}= \frac{2xy\sqrt{17}}{3x+3y+\sqrt{9x^2-16xy+9y^2}},$$
and the radius of the circumcircle is 
$$R= \frac{6y\cdot 2\sqrt{9x^2-16xy+9y^2} \cdot 6x}{2\cdot 4xy\sqrt{17}}=\frac{9\sqrt{9x^2-16xy+9y^2}}{\sqrt{17}}.$$
\end{proof}

\section{Conclusion }
To view Primitive pythagorean triples as coordinates of points in 3-dimensional space opens interesting posibilities of exploring the properties of Pythagorean triples. We intend to continue our study of PPT and bring more of the related results. Namely we would like to focus on describing the metric properties of descendants of PPT in Berggren tree.

\begin{flushright}
	Corresponding author: Lucia JANI\v{C}KOV\'A
	
	\textit{Institute of Mathematics, P.J. \v Saf\'arik University in Ko\v sice, Slovakia}
	
	\textit{Jesenn\'a 5, 041 54 Ko\v sice, Slovakia}
	
	e-mail: \textit{lucia.janickova@upjs.sk}

	Evelin CS\'OK\'ASI
	
	\textit{Institute of Mathematics, P.J. \v Saf\'arik University in Ko\v sice, Slovakia}
	
	\textit{Jesenn\'a 5, 041 54 Ko\v sice, Slovakia}
	
	e-mail: \textit{csokasi.evelinke@gmail.com}
	
\end{flushright}

\end{document}